\documentstyle[auxdefs,diagram,12pt]{article}

\begin{document}

\input amssym.def

\input amssym

\newcommand{\ba}{\mbox{$\Bbb A$}}
\newcommand{\bz}{\mbox{$\Bbb Z$}}
\newcommand{\bq}{\mbox{$\Bbb Q$}}
\newcommand{\bc}{\mbox{$\Bbb C$}}
\newcommand{\bh}{\mbox{$\Bbb H$}}
\newcommand{\br}{\mbox{$\Bbb R$}}
\newcommand{\bp}{\mbox{$\Bbb P$}}
\newcommand{\bres}{\mbox{$\Bbb F$}}

\newcommand{\ra}{\mbox{$\rightarrow$}}
\newcommand{\lra}{\mbox{$\longrightarrow$}}
\newcommand{\proof}{\noindent {\it Proof: }}
\newcommand{\remark}{\noindent {\bf Remark: }}

\newcommand{\cA}{\mbox{${\cal A}$}}
\newcommand{\cB}{\mbox{${\cal B}$}}
\newcommand{\cC}{\mbox{${\cal C}$}}
\newcommand{\cD}{\mbox{${\cal D}$}}
\newcommand{\cE}{\mbox{${\cal E}$}}
\newcommand{\cF}{\mbox{${\cal F}$}}
\newcommand{\cG}{\mbox{${\cal G}$}}
\newcommand{\cH}{\mbox{${\cal H}$}}
\newcommand{\cI}{\mbox{${\cal I}$}}
\newcommand{\cJ}{\mbox{${\cal J}$}}
\newcommand{\cK}{\mbox{${\cal K}$}}
\newcommand{\cL}{\mbox{${\cal L}$}}
\newcommand{\cM}{\mbox{${\cal M}$}}
\newcommand{\cN}{\mbox{${\cal N}$}}
\newcommand{\cO}{\mbox{${\cal O}$}}
\newcommand{\cP}{\mbox{${\cal P}$}}
\newcommand{\cQ}{\mbox{${\cal Q}$}}
\newcommand{\cR}{\mbox{${\cal R}$}}
\newcommand{\cS}{\mbox{${\cal S}$}}
\newcommand{\cT}{\mbox{${\cal T}$}}
\newcommand{\cU}{\mbox{${\cal U}$}}
\newcommand{\cV}{\mbox{${\cal V}$}}
\newcommand{\cW}{\mbox{${\cal W}$}}
\newcommand{\cX}{\mbox{${\cal X}$}}
\newcommand{\cY}{\mbox{${\cal Y}$}}
\newcommand{\cZ}{\mbox{${\cal Z}$}}

\setlength{\unitlength}{1cm} \thicklines \newcommand{\ci}{\circle*{.13}}
\renewcommand{\mp}{\multiput}
\newcommand{\num}[1]{{\raisebox{-.5\unitlength}{\makebox(0,0)[b]{${#1}$}}}}
\newcommand{\mynum}[1]{{\raisebox{-.1\unitlength}{\makebox(-.4,0)[r]{$Y_{#1}$}}}}

\newcommand{\da}{\downarrow}
\newcommand{\el}{\ell}

\newcommand{\aut}{automorphism}
\newcommand{\clg}{compact Lie group}
\newcommand{\cpinf}{\mbox{${\Bbb P}^\infty$}}
\newcommand{\defr}{deformation retract}
\newcommand{\diffeo}{diffeomorphism}
\newcommand{\ems}{Eilenberg-Maclane}
\newcommand{\fg}{finitely-generated} 
\newcommand{\fd}{finite-dimensional}

\newcommand{\hep}{homotopy extension property}
\newcommand{\homeo}{homeomorphism}
\newcommand{\ho}{homomorphism}
\newcommand{\homeq}{homotopy equivalence}
\newcommand{\hty}{homotopy type}
\newcommand{\iso}{isomorphism}
\newcommand{\les}{long exact sequence}
\newcommand{\lfp}{Lefschetz fixed-point theorem}
\newcommand{\mvs}{Mayer-Vietoris sequence}
\newcommand{\ses}{short exact sequence}
\newcommand{\nbhd}{neighborhood}
\newcommand{\pid}{principal ideal domain}
\newcommand{\sss}{spectral sequence}
\newcommand{\we}{weak equivalence}

\newcommand{\wrt}{with respect to}

\newcommand{\rn}{\mbox{$\br ^n$}}
\newcommand{\cn}{\mbox{$\bc ^n$}}

\newcommand{\gc}{\mbox{$G_{\bc}$}}

\newcommand{\phiplus}{\mbox{$\Phi ^+$}}
\newcommand{\phiminus}{\mbox{$\Phi ^-$}}
\newcommand{\phii}{\mbox{$\Phi _I$}}
\newcommand{\phij}{\mbox{$\Phi _J$}}
\newcommand{\phiiplus}{\mbox{$\Phi _I ^+$}}
\newcommand{\phiiminus}{\mbox{$\Phi _I ^-$}}
\newcommand{\phijplus}{\mbox{$\Phi _J ^+$}}
\newcommand{\phijminus}{\mbox{$\Phi _J ^-$}}
\newcommand{\phiirad}{\mbox{$\Phi _I ^{rad}$}}
\newcommand{\phijrad}{\mbox{$\Phi _J ^{rad}$}}
\newcommand{\phiipar}{\mbox{$\Phi _I ^{par}$}}
\newcommand{\phijpar}{\mbox{$\Phi _J ^{par}$}}
\newcommand{\ujrad}{\mbox{$U _J ^{rad}$}}
\newcommand{\uirad}{\mbox{$U _I ^{rad}$}}

\newcommand{\lieuirad}{u_I ^{rad}}
\newcommand{\lieujrad}{u_J ^{rad}}
\newcommand{\lieualpha}{u_{\alpha}}
\newcommand{\lieqw}{q_w}
\newcommand{\liepi}{p_I}
\newcommand{\lieuw}{u_w}
\newcommand{\lieuwp}{u_w ^\prime}
\newcommand{\liegc}{\frak{g}_{\bc}}
\newcommand{\liehc}{\frak{h}_{\bc}} 
\newcommand{\gctil}{\mbox{$\tilde{G} _{\bc}$}}
\newcommand{\gctilp}{\mbox{$\tilde{G} _{\bc}/P$}}
\newcommand{\gctilad}{\mbox{$\tilde{G} _{\bc}^{ad}$}}
\newcommand{\gctiladp}{\mbox{$\tilde{G} _{\bc} ^{ad}/P^{ad}$}}
\newcommand{\btil}{\mbox{$\tilde{B}$}}
\newcommand{\btilminus}{\mbox{$\tilde{B}^-$}}  
\newcommand{\dtil}{\mbox{$\tilde{D}$}}
\newcommand{\util}{\mbox{$\tilde{U}$}}
\newcommand{\utilminus}{\mbox{$\tilde{U} ^-$}}
\newcommand{\utilminussig}{\mbox{$\tilde{U} ^- _\sigma$}}
\newcommand{\wtil}{\mbox{$\tilde{W}$}}
\newcommand{\wtils}{\mbox{$\tilde{W}^S$}}
\newcommand{\wtili}{\mbox{$\tilde{W}^I$}}   
\newcommand{\wtilil}{\mbox{$\tilde{W}_I$}}   
\newcommand{\ntil}{\mbox{$\tilde{N}$}}
\newcommand{\stil}{\mbox{$\tilde{S}$}}
\newcommand{\liegctil}{\mbox{$\tilde{\frak{g}} _{\bc}$}}
\newcommand{\phitil}{\mbox{$\tilde{\Phi}$}}
\newcommand{\phitilplus}{\mbox{$\tilde{\Phi} ^+$}}
\newcommand{\phitilminus}{\mbox{$\tilde{\Phi} ^-$}}
\newcommand{\coroot}{\mbox{$Q^\vee$}}
\newcommand{\antidom}{\mbox{$Q^\vee _-$}}
\newcommand{\chamber}{\overline{\cC}}
\newcommand{\affg}{\mbox{$\cL _G$}}   
\newcommand{\tchat}{\mbox{$\hat{T} _{\bc}$}}
\newcommand{\that}{\mbox{$\hat{T}$}}  

\newcommand{\afv}{affine flag variety}
\newcommand{\rvar}{Richardson variety}
\newcommand{\rvars}{Richardson varieties}  
\newcommand{\rfilt}{Richardson filtration} 
\newcommand{\bvar}{Birkhoff variety}
\newcommand{\bvars}{Birkhoff varieties}
\newcommand{\cpo}{closed parabolic orbit}
\newcommand{\clr}{coroot lattice representative}  
\newcommand{\ef}{equivariant formality}
\newcommand{\svar}{Schubert variety}
\newcommand{\svars}{Schubert varieties}
\newcommand{\pd}{Poincar\' e duality}
\newcommand{\pds}{Poincar\' e duality space}
\newcommand{\ppoly}{Poincar\' e polynomial}
\newcommand{\pser}{Poincar\' e series}
\newcommand{\psg}{one-parameter subgroup}
\newcommand{\fnk}{\mbox{$X_{n,k}$}}
\newcommand{\fnkp}{\mbox{$X_{n,k} ^\flat$}}
\newcommand{\ppc}{positive pair condition}
\newcommand{\npc}{negative pair condition}
\newcommand{\osc}{opposite sign condition}
\newcommand{\und}{\underline{0}}
\newcommand{\nope}{$\lambda \notin \cQ ^\vee$}
\newcommand{\onemin}{1- \alpha _0 (\lambda)} 
\newcommand{\owt}{(otherwise $\lambda$ is overweight)} 
\newcommand{\adp}{admissible palindromic}
\newcommand{\fork}{(otherwise $\lambda$ covers a fork)}

\newcommand{\gsm}{graph-splitting move} 
\newcommand{\aslam}{\alpha _s
(\lambda)} 
\newcommand{\azlam}{\alpha _0 (\lambda)}
\newcommand{\atlam}{\alpha _t (\lambda)} 
\newcommand{\aulam}{\alpha _u
(\lambda)} 
\newcommand{\astlam}{\alpha _s (\lambda) + \alpha _t
(\lambda)} 
\newcommand{\alam}{\alpha (\lambda)} 
\newcommand{\blam}{\beta
(\lambda)} 
\newcommand{\aplam}{\alpha ^\prime (\lambda)}
\newcommand{\bpp}{$\beta$-positive pair}
\newcommand{\bnp}{$\beta$-negative pair}
\newcommand{\lam}{\mbox{$\lambda$}}
\newcommand{\circular}{spiral}

\newcommand{\xlam}{\mbox{$X_\lambda$}}
\newcommand{\zlam}{\mbox{$Z_\lambda$}}
\newcommand{\elam}{\mbox{$e_\lambda$}}
\newcommand{\slam}{\mbox{$S_\lambda$}}
\newcommand{\openlam}{\mbox{$\cU _\lambda$}}
\newcommand{\openslam}{\mbox{$\cS _\lambda$}}
\newcommand{\nlam}{\mbox{$N_\lambda$}}
\newcommand{\xmu}{\mbox{$X_\mu$}}
\newcommand{\zmu}{\mbox{$Z_\mu$}}
\newcommand{\emu}{\mbox{$e_\mu$}}
\newcommand{\smu}{\mbox{$S_\mu$}}
\newcommand{\openmu}{\mbox{$\cU _\mu$}}

\newcommand{\vlamu}{\mbox{$V_{\lambda ,\mu}$}}
\newcommand{\klamu}{\mbox{$K_{\lambda ,\mu}$}}
\newcommand{\vlamuplus}{\mbox{$V_{\lambda ,\mu} ^+$}}
\newcommand{\vlam}{\mbox{$\cE _\lambda$}}
\newcommand{\vi}{\mbox{$\cE _{\cal I}$}}
\newcommand{\zi}{\mbox{$Z _{\cal I}$}}
\newcommand{\zip}{\mbox{$Z _{{\cal I} ^\prime}$}}
\newcommand{\xj}{\mbox{$X _{\cal J}$}}
\newcommand{\sj}{\mbox{${\cal S} _{\cal J}$}}   

\newcommand{\utilam}{\mbox{$\tilde{U} _\lambda$}}
\newcommand{\utilaminus}{\mbox{$\tilde{U} _\lambda ^-$}}
\newcommand{\utilamp}{\mbox{$\tilde{U} _\lambda ^\prime$}} 
 \newcommand{\utilampminus}{\mbox{$\tilde{U} _\lambda ^{\prime -}$}} 
\newcommand{\lieutilam}{\mbox{$\frak{u} _\lambda ^\prime$}}  
\newcommand{\lieutilamp}{\mbox{$\tilde{U} _\lambda$}}  
\newcommand{\ulamp}{\mbox{$\tilde{U} _\lambda$}}

\newtheorem{theorem}{Theorem}[section]

\newtheorem{proposition}[theorem]{Proposition}

\newtheorem{lemma}[theorem]{Lemma}

\newtheorem{conjecture}[theorem]{Conjecture}

\newtheorem{example}[theorem]{Example}

\newtheorem{corollary}[theorem]{Corollary}

\newtheorem{definition}[theorem]{Definition}

\title{The topology of Birkhoff varieties}

\author{Luke Gutzwiller and Stephen A. Mitchell}

\maketitle

\section{Introduction}

Let $\cF =\tilde{G}/P_I$ be an affine flag variety. Here $G$ is a
simply-connected complex algebraic group with simple Lie algebra,
$\tilde{G} =G(\bc [z,z^{-1}])$ is the corresponding affine group, and
$P_I$ is the parabolic subgroup associated to a subset $I$ of the set of
Coxeter generators \stil\ of the affine Weyl group \wtil . Then
\cF\ has two dual stratifications: the Schubert or Bruhat cell
decomposition

$$\cF =\coprod _{\lambda \in \tilde{W} /\tilde{W}_I} e_\lambda ,$$

\noindent and the Birkhoff stratification

$$\cF =\coprod _{\lambda \in \tilde{W} /\tilde{W}_I} S_\lambda .$$ 
 
The Schubert cells \elam\ are the orbits of the Iwahori subgroup \btil ,
while the Birkhoff strata \slam\ are the orbits of the opposite Iwahori
subgroup \btilminus . The cells and the strata are dual in the sense
that $\slam \cap \elam =\{\lam \}$, and the intersection is
transverse. The closure of \elam\ is the affine Schubert variety \xlam
. It has dimension $\el ^I (\lambda)$, where $\el ^I$ is the minimal
length occuring in the coset $\lam \wtilil$, and its cells are indexed
by the lower order ideal generated by \lam\ in the Bruhat order on
$\wtil /\wtilil$. Dually, the closure of \slam\ is the Birkhoff variety
\zlam . It is an infinite-dimensional irreducible ind-variety with
codimension $\el ^I (\lambda)$. Its Birkhoff strata are indexed by the
upper order ideal generated by \lam .

Thus the Birkhoff varieties may be viewed as analogous to the dual
\svars\ from the classical setting, in which the role of \cF\ is
played by a finite-dimensional flag variety. More generally, let $\cI$
denote an upper order ideal in the Bruhat poset $\wtil /\wtilil$. Then
$\zi =\cup _{\lambda \in {\cal I}} \slam$ is a finite union of \bvars . 
Our main theorem shows that in one respect, the classical and affine
cases differ dramatically.

\begin{theorem} \label{main} \marginpar{}

Let \zi\ be a finite union of \bvars\ in the \afv\ \cF . Then \zi\ is a
\defr\ of \cF . In particular, the inclusion $\zi \subset \cF$ induces
\iso s on ordinary and equivariant cohomology, with any coefficients.

\end{theorem}
 
The proof has two main ingredients. The first is the existence of a sort of
``algebraic tubular neighborhood'' of \zi . Let $\vi =\cup _{\lambda
\in {\cal I}} \elam $. Then \vi\ is a Zariski open \nbhd\ of \zi
. Similarly, let \cJ\ be a proper lower order ideal in $\wtil /\wtilil$,
let $\xj =\cup _{\lambda \in {\cal J}} \elam$, and let $\sj =\cup
_{\lambda \in {\cal J}} \slam$. Then $\xj$ is a finite union of \svars ,
and $\cS _\lambda$ is a Zariski open \nbhd\ of \xj . Then the following
theorem holds for both affine and classical flag varieties.

\begin{theorem} \label{} \marginpar{}

a) \zi\ is a \defr\ of \vi . 

b) \xj\ is a \defr\ of \sj .  

\end{theorem}

Versions of part (b) appear to be known (see for example
the special case discussed in \cite{cls}), but we are not aware of a
proof or even a full statement of this theorem in the literature.

The second ingredient depends on the infinite-dimensionality of the
Birkhoff strata, and has no analog in the classical case. 

\begin{lemma} \label{} \marginpar{}

The punctured Birkhoff stratum $\slam -\{\lam \}$ is contractible. 

\end{lemma}

Given these two ingredients, the main theorem follows by a formal
downward induction over the Birkhoff filtration, using 
Whitehead's theorem at the inductive step. 

\bigskip

\noindent {\it Organization of the paper:} In \S 2 we summarize some
basic notation, and introduce a well-known $\bc ^\times$-action or
complex flow on \cF\ that will be used to construct our
deformations.  In \S 3 we study complex flows on ind-spaces and
ind-varieties. The main result is a general criterion for deforming an
ind-space into an invariant ind-subspace using a flow
(Theorem~\ref{flow}). In \S 4 we study the structure of the affine
analog \util\ of a maximal unipotent subgroup, and its opposite
\utilminus . The main application is to show that punctured Birkhoff
strata are contractible (Lemma~\ref{puncture}). In \S 5 we construct our
algebraic tubular \nbhd s (Theorem~\ref{tube}). 

In \S 6 we prove the main theorem. We also compute the homology of the
pairs $(\vi , \vi -\zi )$ and $(\sj , \sj -\xj)$. These pairs can be
viewed as algebraic normal Thom spaces of the ind-subvarieties \zi , \xj\ in
\cF . Finally, we make some remarks on torus-equivariant cohomology
$H^* _{\hat{T}} \zi$. In particular, we prove one half of a
Goresky-Kottwitz-MacPherson theorem (Proposition~\ref{equi}). 

\bigskip

\noindent {\it Acknowledgements:} We would like to thank Aravind Asok,
Sara Billey, Megumi Harada, and Shrawan Kumar for helpful conversations.  
We would also like to think the referee for some valuable
suggestions and questions.

\section{Preliminaries}

We use the following conventions throughout this paper: 

\bigskip

All (co)homology groups are singular (co)homology groups with integer
coefficients, unless otherwise specified. 

Varieties over \bc\ are given the classical Hausdorff topology inherited
from $\bc ^n$ or $\bp ^n$, which we call the {\it complex}
topology. When the Zariski topology is used, it will be indicated
explicitly. Likewise, ind-varieties have both a complex and a Zariski
direct limit topology.

The term {\it deformation retract} means what some authors call {\it
  strong deformation retract}; i.e., the deformation fixes the subspace
  in question pointwise. 

\subsection{Notation}

\noindent {\it The group} $G$. Let $G$ be a simply-connected complex
algebraic group with simple Lie algebra, with maximal torus $T_{\Bbb
C}$, Weyl group $W$, $S \subset W$ the simple reflections, root system
$\Phi$, and simple roots $\alpha _s, s\in S$. Let \coroot\ denote the
coroot lattice. Let $B$ denote a Borel subgroup containing $T_{\bc}$,
and $U \subset B$ the unipotent radical. Let $B^-, U^-$ denote the
opposite Borel and unipotent subgroups. We write $\frak{g}, \frak{u} $,
and so on for the Lie algebras.

\bigskip

\noindent {\it Affine groups.}  Let $\tilde{G} =G(\bc [z,z^{-1}])$; this
is the group of regular maps $\bc^\times \lra G$. Similarly $P=G(\bc
[z])$ is the group of regular maps $\bc \lra G$. We have subgroups $P
\supset \tilde{B} \supset \util \supset P^{(1)}$ defined as follows: The
Iwahori subgroup is $\btil =\{f \in P: f(0)\in B^-\}$; similarly $\util
=\{f \in P: f(0) \in U^-\}$. Set $P^{(1)} =\{f \in P: f(0)=1\}$. Let
$P^-=G(\bc [z^{-1}])$ denote the group of regular maps $\bp ^1 -\{0\}
\lra G$. Analogs of the subgroups of $P$ are defined in the evident way;
e.g. $\tilde{B}^-=\{f \in P: f(\infty)\in B\}$, etc. Associated Lie
algebras are written $\tilde{\frak{g}} =\frak{g} \otimes \bc
[z,z^{-1}]$, and so on.

The group $\tilde{G}$ is an affine ind-group. Explicitly, in the case $G
=SL_n \bc$ we let $F_m SL_n \bc [z,z^{-1}]$ denote the subset of
matrices $A$ such that $A_{ij} =\sum _{k=-m} ^m a_{ijk} z^k$. This
defines a filtration $F_m$ by affine varieties that yields the affine
ind-group structure. In the general case we choose a faithful
representation $G \subset SL_n \bc$ and set $F_m \tilde{G} =\tilde{G} \cap
F_m SL_n \bc [z,z^{-1}]$. It is easy to see that the affine ind-group
structure obtained is independent of the choice of representation. For a
more general Kac-Moody approach, see \cite{kumar}, \S 7.3.

\bigskip

\noindent {\it Affine Weyl group.}  Let \wtil\ denote the affine Weyl
group, with Coxeter generators $\stil =S \cup \{s_0\}$. The affine root
system is $\phitil =\bz \times \Phi$. As simple system of positive roots
we take $\{ (0, -\alpha _s): s \in S\} \cup \{ (1,\alpha _0)\}$, where
$\alpha _0$ is the highest root. If $\theta =(n, \alpha)$, let $r_\theta
=r_{n, \alpha}$ denote the affine reflection associated to
$(n,\alpha)$. 

The affine roots occur as weights of the extended torus $\tchat =\bc
^\times \times T_{\bc}$ acting on $\tilde{\frak{g}}$. Here the extra
factor $\bc ^\times$ is acting by loop rotation. Thus $\phitil$ is
actually the set of so-called ``real'' roots; we will also need the
``imaginary'' roots $(n,0)$, $n \in \bz -\{0\}$, which are the weights
of the \tchat\ action on $\frak{t}_{\bc} \otimes \bc \cdot z^n \subset
\tilde{\frak{g}}$.

\bigskip

\noindent $\wtil /\wtilil$ {\it and Bruhat order.} Let \wtili\ denote
the set of minimal length representatives for the cosets $\wtil
/\wtilil$. For any $\sigma
\in \wtil$, let $\el ^I (\sigma)$ denote the $I$-length of $\sigma$;
that is, the length of the minimal coset representative in $\sigma
W_I$. 
Let $\cI ^\lambda$ (resp. $\cJ _\lambda$) denote the upper order
ideal (resp. lower order ideal) generated by \lam\ in the Bruhat order
$\leq$ on $\wtil /W_I$.  We write $\lambda \da \mu$ when $\mu <\lambda$
and the $I$-lengths differ by 1.

\bigskip

\noindent {\it Parabolic subgroups.} Let $P_I \subset \tilde{G}$ denote
the parabolic subgroup generated by \btil\ and $I$. Then $P_I$ is the
semi-direct product of a normal subgroup $\tilde{U} _I$ and a \fd\
subgroup $L_I$. Here $\tilde{U} _I \subset \tilde{U}$ plays the role of
unipotent radical, and $L_I$ is the Levi factor. Similarly, the opposite
parabolic $P_I ^-$ generated by $\tilde{B} ^-$
and $I$ is the semi-direct product of $L_I$ and a normal subgroup
$\tilde{U} ^- _I$.

\bigskip

\noindent {\it Affine flag varieties.} An affine flag variety is
homogeneous space of the form $\cF=\tilde{G} /P_I$. It has a canonical
structure of projective ind-variety (\cite{kumar}, 13.2.13-18,
\cite{ps}). Set $\cU _0 =\utilminus P/P$ and $\openlam =\lambda \cU _0$
for $\lam \in \wtil /\wtilil$ (note this is well-defined). Then the natural
map $U_I^- \lra \cU _0$ is an \iso\ of ind-varieties, and the
$\openlam$'s form a Zariski open cover of \cF . The {\it Birkhoff
strata} \slam\ are the orbits of \btilminus\ on \cF .

\bigskip

\noindent {\it Schubert and Birkhoff varieties.} It is easy to see that
any infinite subset of \wtili\ is cofinal for the Bruhat order
(cf. \cite{bm}, Proposition 7.1). Hence any proper lower order ideal
\cJ\ is finite, and $\xj =\cup _{\sigma \in {\cal J}} e_\sigma$ is a
finite union of \svars . If $\cJ =\cJ _\lambda$, this is just the \svar\
\xlam.  If \cI\ is any non-empty upper order ideal, then $\zi =\cup
_{\sigma \in {\cal I}} S_\sigma$ is a finite union of Birkhoff
varieties. When $\cI =\cI ^\lambda$, this is just the \bvar\ \zlam .

Define $\utilam =\util \cap \lambda \tilde{U}_I ^- \lambda ^{-1}$, 
$\utilamp =\util \cap \lambda P_I \lambda ^{-1}$,
$\utilaminus =\tilde{U}^- \cap \lambda \tilde{U}_I ^- \lambda ^{-1}$,
$\tilde{U}_\lambda ^{\prime -}=\tilde{U}^- \cap \lambda P_I \lambda
^{-1}$. Thus \utilamp\ and $\tilde{U}_\lambda ^{\prime -}$ are the
isotropy groups of the \util\ and \utilminus\ actions on $\lambda P_I
/P_I$, while the group action defines \iso s $\utilam \cong e_\lambda$
and $\utilaminus \cong \slam$.

\subsection{The extended torus action and the flow}

 Let $\tchat$ denote the extended torus $\bc ^\times \times
T_{\bc}$. Then \tchat\ acts on $\tilde{G}$: The constant torus valued
loops $T_{\bc}$ act by conjugation, while the extra factor $\bc ^\times$
acts by loop rotation. The action preserves parabolic subgroups and
induces an algebraic group action $\tchat \times \cF \lra \cF$, with
fixed point set $\wtil /\wtilil $. The action also preserves \btilminus
, Schubert cells, Birkhoff strata, etc.  The action of \tchat\ on a
Schubert cell \elam\ is isomorphic to a linear action, with weights
precisely the set of roots $\tilde{\Phi} _\lambda$, each occuring with
multiplicity one. In particular, the weights are positive.
 
Now consider the action of the torus $\hat{T}=\bc ^\times \times
T_{\bc}$ on \cF . One can always find a rank one subtorus $\phi :
\bc ^\times \lra \hat{T}$ such that the induced $\bc ^\times$ action has
the following properties:

\bigskip                                                                       

(i) The fixed-point set is still $\wtil/W$;

(ii) If $x \in \elam$, then $lim _{t \ra 0} t \cdot x=\lambda$;

(iii) If $x \in \slam$, then $lim _{t \ra \infty} t \cdot x=\lambda$.          

\bigskip                                                                        
To see this, identify $Hom \, (\bc ^\times , \hat{T})$ with $\bz \times        
\coroot$ and write $\phi =(k, \gamma)$. We then have: 

\begin{proposition} \label{} \marginpar{}                                       
Suppose that (a) For all $\alpha \in \phiplus$, $\alpha (\gamma)<0$, and
(b) $k> max _{\alpha \in \Phi} |\alpha (\gamma)|$. Then $\phi
=(k,\gamma)$ has properties (i)-(iii) above.

In particular, (i)-(iii) hold when $\gamma =-\sum _{\alpha \in \phiplus}
\alpha ^\vee$ and $k=2h-1$, where $h$ is the Coxeter number.
\end{proposition}                                                               
\proof Assumptions (a) and (b) ensure that $\bc ^\times $ acts on each
cell $e_\lambda$ with positive weights, yielding (i) and (ii). Now
suppose $x\in \utilminus \lambda P_I/P_I$. Since $\utilminus$ is generated
by the root subgroups $U_{n,\alpha}$ with $(n,\alpha) \in \phitilminus$
\cite{kp}, and $\bc ^\times$ acts on these with negative weights, it
follows that $lim _{t\ra \infty} t\cdot x=\lambda$, proving (iii).  For
the last assertion of the proposition, let $\rho ^\vee =\omega _1 ^\vee
+...+\omega _r ^\vee$, where the $\omega _i ^\vee$'s are the fundamental
coweights and $r$ is the rank of $G$. Let $\alpha _0 =\sum _{i=1}^r m_i
\alpha _i$, where the $\alpha _i$'s are the simple positive roots and
$\alpha _0$ is the highest root as usual. Then $\gamma =-2\rho ^\vee$,
verifying (a), while the max occurring in (b) is $\alpha _0 (2 \rho
^\vee) =2 \sum m_i =2h-2$.

\bigskip

Fix $\gamma$, $k$ as in the Proposition. We refer to the resulting $\bc
^\times$ action as the {\it complex flow}. 

\section{Ind-spaces, ind-varieties and $\bc^\times$-actions}

\subsection{Ind-spaces and ind-varieties}

An {\it ind-space} is a set $X$ equipped with a filtration $X_1 \subset
X_2 \subset ...$ such that $X=\cup X_n$, each $X_n$ is a topological
space, and $X_n$ is closed in $X_{n+1}$. We give $X$ the direct limit
topology: A subset of $X$ is closed if and only if its intersection with
each $X_n$ is closed. A morphism of ind-spaces is a map $f: X\lra Y$
such that for every $n$ there exists $m$ with $f(X_n) \subset Y_m$ and
$f: X_n \lra Y_m$ continuous. In particular, $f$ is continuous.  Two
ind-space structures on the same space $X$ are {\it commensurate} if the
identity map is an \iso\ between them. Any subspace $A$ of an ind-space
$X$ is an ind-space with $A_n =A \cap X_n$. Given any space $X$, we can
form the {\it constant ind-space} with $X_n=X$ for all $n$. This embeds
the category of spaces as a full subcategory of the category of
ind-spaces.

An {\it ind-variety} is defined similarly, with the requirement that
  each $X_n$ is a complex algebraic variety and a closed subvariety of
  $X_{n+1}$. See \cite{kumar} for a brief introduction to
  ind-varieties. Every ind-variety is an ind-space in the Zariski and
  complex topologies. An ind-variety is {\it irreducible} if it is
  irreducible as a topological space in the Zariski topology. If each
  filtrant $X_n$ is irreducible, then so is $X$. Conversely, if $X$ is
  irreducible then it admits a commensurate filtration $Y_n$ with each
  $Y_n$ irreducible. In fact any filtration with $Y_n$ an irreducible
  component of $X_n$ is a commensurate filtration. 

If $H$ is an algebraic group, an {\it ind-H-variety} is an ind-variety $X$
equipped with compatible algebraic $H$-actions on each $X_n$. If $H$ is
connected and $V$ is any $H$-variety, each irreducible component of $V$
is invariant under the action. It follows that if $H$ is connected, then
any ind-H-variety has a commensurate filtration by irreducible
$H$-invariant varieties. 

A {\it group ind-variety}, or simply ``ind-group'', is a group object
$\Gamma$ in the category of ind-varieties. Note that the filtrants
$\Gamma _n$ are not assumed to be subgroups. A connected ind-group is
irreducible (\cite{kumar}, Lemma 4.2.5). 

\subsection{Ind-CW-complexes}

An {\it ind-CW-complex} is an ind-space such that each $X_n$ admits a
CW-structure having $X_{n-1}$ as a subcomplex. We do not assume these
structures are compatible as $n$ varies, and indeed $X$ itself need not
admit any CW-structure (see the example below). An {\it ind-CW pair} is
a pair of ind-spaces $(X,A)$ such that each $X_n$ admits a CW-structure
such that $X_{n-1}$ and $A_n$ are subcomplexes. In the complex topology
an ind-variety is also an ind-CW-complex, by Hironaka's theorem
\cite{hironaka}. However, an ind-variety need not admit any
CW-structure.

\bigskip

\noindent {\it Example:} For $n \geq 1$ let $H_n$ denote the hyperplane
$x=1/n$ in $\bc ^2$. Let $X_n$ denote the union of the coordinate axes
and $H_1,...,H_n$. Let $X=\cup _n X_n$, with its evident ind-variety
structure. Then $X$ does not admit a CW-structure. To see this, suppose
given a CW-structure on $X$, and let $p_n =(1/n,0)$. Then no $p_n$ lies
in a 2-cell, since $X-\{p_n\}$ is disconnected. Furthermore, only
finitely many $p_n$'s can be vertices, since $p_n \lra (0,0)$ in the
direct limit topology, and the vertex set of a CW-complex has no limit
points. More generally, a subset of a CW-complex whose intersection with
each cell is finite has no limit points. Thus all but finitely many
$p_n$'s must lie in a single 1-cell $e^1$. Let $\phi : (0,1) \lra e^1$
be a \homeo . Then for some $n$ we have $a_{-1} < a_0 <a_1 \in (0,1)$
with $\phi (a_i) =p_{n+i}$. This forces $e^1 \cap H_n =p_n$, since if
the path $\phi$ ever enters $H_n -\{p_n\}$ then it must also exit through
$p_n$, contradicting the injectivity of $\phi$. But if $e^1 \cap H_n
=p_n$, then no vertex of $H_n -\{p_n\}$ can be connected by a 1-cell to
$p_n$. Hence the 1-skeleton of $X$ is disconnected, contradicting the
connectedness of $X$. 

\bigskip

As the following two results illustrate, however, for many purposes
ind-CW-complexes are just as good as CW-complexes. 

\begin{lemma} \label{} \marginpar{}

Let $(X,A)$ be an ind-CW-pair. Then $(X,A)$ has the \hep . 

\end{lemma}

\proof See \cite{hatcher}, Chapter 0 for a discussion of the \hep . Any
CW-pair has the \hep ; the lemma follows immediately by an induction
argument, using the CW-pair $(X_{n+1} , X_n \cup A_{n+1})$ at the
inductive step. 

\bigskip

A {\it CW-space} is a space with the homotopy-type of a CW-complex. 

\begin{proposition} \label{} \marginpar{}

Let $X$ be an ind-CW-complex. Then $X$ is a CW-space. 

\end{proposition}

\proof For any space $Y$, there is a CW-approximation $\eta _Y : W(Y)
\lra Y$; that is, a CW-complex $W(Y)$ and a \we\ $\eta _Y$ (see
\cite{hatcher}, Chapter 4). In fact one can make $W$ a functor and
$\eta$ a natural transformation from $W$ to $Id$, by taking $W(Y)$ to be
the geometric realization of the singular complex of $Y$. Hence there is
a functorial CW-approximation $\eta _X : W(X) \lra X$ that is filtered
by CW-approximations $W(X_n) \lra X_n$, with $W(X_n)$ a subcomplex of
$W(X_{n+1})$. By Whitehead's theorem (\cite{hatcher}, Theorem 4.5), each
$W(X_n) \lra X_n$ is a homotopy equivalence. Since each of the pairs
$(X_{n+1}, X_n)$ and $(W(X_{n+1}), W(X_n))$ has the \hep , it follows by
a standard argument that the direct limit map $W(X) \lra X$ is also a
homotopy equivalence (see \cite{hatcher}, Proposition 4G1 and the
paragraph following its proof).

\bigskip

Thus Whitehead's theorem applies to ind-CW-complexes. In particular, we
have: 

\begin{corollary} \label{wedefr} \marginpar{}

Let $(X,A)$ be an ind-CW-pair, and suppose the inclusion $i: A \subset
X$ is a \we . Then $A$ is a \defr\ of $X$. In particular $A$ is a \defr\
of $X$ if $A,X$ are simply-connected and $H_* i$ is an \iso .

\end{corollary}

\proof By Whitehead's theorem, $i$ is a homotopy equivalence. Since
$(X,A)$ has the \hep , the first conclusion follows from \cite{hatcher},
Corollary 0.20. If $A,X$ are simply-connected and $H_* i$ is an \iso ,
then $i$ is automatically a \we\ (\cite{hatcher}, Corollary 4.33). 

\subsection{$\bc ^\times$ actions}

Let $X$ be an ind-space with $\bc ^\times$ action such that each
  filtration $X_n$ is invariant under the action. We also call this a
  {\it complex flow}. If $\cW \subset X$ is open and $C$ is any subset
  of $X$, we say that $C$ {\it flows to} \cW\ {\it at zero} if for every
  $n$ there is an $s>0$ such that for all $|t|\leq s$ we have $t \cdot
  C_n \subset \cW _n$. We say that $C$ {\it flows to} \cW\ {\it at} $\infty$
  if for every $n$ there is an $s>0$ such that for all $|t|\geq s$ we
  have $t \cdot C_n \subset \cW _n$.

A closed $\bc ^\times$-invariant ind-subspace $A$ is {\it strongly
  attractive at zero} (resp. {\it strongly attractive at $\infty$}) if
  for every \nbhd\ \cW\ of $A$ and $x \in X$, there is a \nbhd\ $U$ of
  $x$ that flows to \cW\ at zero (resp. at $\infty$).  Since the
  conditions ``attractive at zero'' and ``attractive at $\infty$'' are
  interchanged under the automorphism $t \lra t^{-1}$ of $\bc ^\times$,
  for the remainder of this section we will consider only the former
  case and call such a subspace {\it strongly attractive}.

\bigskip

\noindent {\it Remark:} Call $A$ {\it weakly attractive} in $X$ if the
above condition merely holds pointwise, i.e., for every \nbhd\ \cW\ of
$A$ and $x \in X$, there is an $s>0$ such for all $|t| \leq s$ we have
$t \cdot x \in \cW$. This is a very weak condition that does not imply
strongly attractive, even if $X$ is a compact constant ind-space and one
adds the requirement that $lim _{t \ra 0} t \cdot x$ exists for all $x$.
For example, take $X=\bp ^1$ with the standard $\bc ^\times$ action
coming from diagonal matrices in $SL_2 \bc$, and take $A$ to consist of
the two fixed points $p_0, p_\infty$. Here we have labelled the points
so that for any $x \notin A$, $t \cdot x \ra p_0$ (resp. $p_\infty$) as
$t \ra 0$ (resp. $\infty$). Then $A$ is weakly attractive but evidently
not strongly attractive (take \cW\ to be the union of disjoint \nbhd s
of $p_0, p_\infty$, and take $x=p_\infty$). One can easily exhibit
similar examples with $A$ connected, for example with $X=\bp ^2$ and
$A=\bp ^1 \bigvee \bp ^1$.

\bigskip

By a {\it regular neighborhood} of a subspace $B$ in a space $Y$, we
mean a \nbhd\ \cW\ such that $B$ is a \defr\ of \cW .

\begin{theorem} \label{flow} \marginpar{}

Let $X$ be a $T_1$ (points are closed) ind-space with $\bc ^\times$
action.  Suppose $A \subset X$ is strongly attractive and each $A_n$ has
a regular \nbhd\ in $X_n$. Then the inclusion $i: A \lra X$ is a \we .
If in addition $(X,A)$ is an ind-CW pair, then $A$ is a \defr\ of $X$.

\end{theorem}

\proof Recall that a \we\ is a map inducing a bijection on
path-components, and an isomorphism on homotopy groups for any choice of
basepoint. We will show that for any compact space $K$, the inclusion
induces a bijection on homotopy classes $i_*: [K, A]
\stackrel{\cong}{\lra} [K,X]$. It is well-known, and easy to prove, that
this implies $i$ is a \we .

Suppose that $X$ is a constant ind-space. Let $f: K \lra X$ be a map,
and let $\cW$ be a regular \nbhd\ of $A$. For each $k \in K$, choose a
\nbhd\ $U_k$ of $k$ and $s_k >0$ such that for all $|t| \leq s_k$, $t
\cdot U_k \subset \cW$. Since $K$ is compact, $f(K)$ is covered by
finitely many such \nbhd s, say $U_{k_1},...,U_{k_n}$. Taking $s=min \,
\{s_{k_1},...,s_{k_n} \}$, we have $s \cdot f (K) \subset \cW$. Since $s
\cdot f$ is homotopic to $f$, composing with the deformation of \cW\
into $A$ shows that $f$ is homotopic to a map $g: K \lra A$. Hence $i_*$
is surjective. Next suppose that $f_0, f_1 : K \lra A$ are maps that
become homotopic in $X$. Applying the preceeding argument to the
homotopy shows that $s \cdot f_0$ is homotopic to $s \cdot f_1$ in $A$,
and hence $f_0$ is homotopic to $f_1$. This shows that $i_*$ is
injective, and hence bijective.

In the general case, we conclude that each inclusion $A_n \subset X_n$
is a \we . Now let $Map \, (-,-)$ denote the set of continuous
maps. Then it is well-known and easy to prove that for any
$T_1$-ind-space $X$ and compact space $K$, the natural map 

$$colim _n Map \, (K, X_n) \lra Map \, (K, X)$$

\noindent is bijective (the $T_1$ hypothesis ensures that every compact
   subset of $X$ lies in some $X_n$). Hence $i_* : [K,A] \lra [K,X]$ is
   a colimit of bijections and so is bijective.

Finally, if $(X, A)$ is an ind-CW-pair then $A$ is a \defr\ of
$X$ by Corollary~\ref{wedefr}.  

\bigskip

\noindent {\it Remark:} Note that the theorem fails miserably if one
only assumes $A$ is weakly attractive in $X$ (see the example in the
previous remark).

\bigskip

When $X$ is an ind-variety, we always assume that the $\bc
^\times$-action is algebraic. The following technical lemma will be need
in the proof of Theorem~\ref{tube}:

\begin{lemma} \label{star} \marginpar{}

Let $f: (X,A) \lra (Y,B)$ be a map of ind-variety pairs with $\bc
^\times$ action. Suppose $f: X \lra Y$ is surjective and satisfies the
following condition:
 
(*) $X$ is a union of ind-subvarieties $Z_\alpha$ such that for each
$\alpha$, the restriction $f |_{Z_\alpha}$ is an \iso\ of
ind-varieties onto a Zariski open ind-subvariety of $Y$. 

Then if $A$ is strongly attractive in $X$, $B$ is strongly attractive
in $Y$.

\end{lemma}

\proof Let $\cW \subset Y$ be a \nbhd\ of $B$, and $y \in Y$. Choose $x
\in f^{-1} y$. Then there is a \nbhd\ $U$ of $x$ that flows to $f^{-1}
\cW$. Moreover, $x \in Z$ for some $Z=Z_\alpha$ as in the theorem. By
Chevalley's theorem each $f(Z_m) \cap Y_n$ is a constructible subset of
$Y_n$, so for fixed $n$ we have $f(Z) \cap Y_n = f(Z_m) \cap Y_n$ for
sufficiently large $m$ (see \cite{kumar}, exercise 7.3.E(2)). Let $V=f(U
\cap Z)$. Then $V$ is complex open, since any \iso\ of varieties is a
\homeo\ in the complex topology. Furthermore, for fixed $n$ and $m>>0$
we have

$$V_n=f(U _m \cap Z) \cap Y_n =f(U \cap Z_m ) \cap Y_n =f(U \cap
Z) \cap Y_n,$$

\noindent where the third equality uses the fact that $f|_Z$ is
injective. Then there is an $s>0$ such that $t \cdot U_m \subset f^{-1}
\cW $ for all $|t| \leq s$, and hence $t \cdot V_n \subset \cW$. Thus
$V$ is a complex open \nbhd\ of $y$ that flows to $\cW$, as required.

\section{\util\ and \utilminus\ as ind-varieties}

In this section we study the structure of \util\ and \utilminus\ as
ind-varieties. In particular, we construct filtrations by weighted
cones. Our main applications are Lemma~\ref{puncture}, showing that
punctured Birkhoff strata are contractible, and
Corollary~\ref{attractpoint}, which will be used in the proof of
Theorem~\ref{tube}. In fact it will suffice to consider \util , for the
following reason: Define $\delta : \tilde{G} \lra \tilde{G} $ by
$(\delta f)(z)=w_0 f(z^{-1})w_0$, where $w_0$ is the longest element of
$W$. Then $\delta$ is an ind-group automorphism exchanging \util\ and
\utilminus . We will leave it to the reader to make the translation from
\util\ to \utilminus ; in particular one must replace positive weights by
negative weights, and limits as $z \ra 0$ by limits as $z \ra \infty$.

\subsection{Weighted cones}

Let $V$ be a \fd\ representation of $\bc ^\times$. By a {\it weighted
  cone} we mean a nonempty, closed $\bc ^\times$-invariant subvariety of
  $V$. We will only be concerned with positively or negatively weighted
  cones. Since the two cases are exchanged under the automorphism of $z
  \mapsto z^{-1}$ of $\bc ^\times$, there is no loss of generality in
  restricting to positively weighted cones. 

\subsubsection{Filtrations by weighted cones}

In this section we show that \util\ is filtered by
  positively  weighted cones. Consider the
  finite-dimensional filtrations $F_m \util \subset \util$,
  defined in \S 2. Recall that $F_m \util$ is an affine variety,
  but not a subgroup. Recall also that \util\ is not unipotent but
  embeds in an inverse limit of unipotent groups. More precisely, it
  embeds in an inverse limit of the form $lim _k U[k]$, where $U[k]$ is
  a maximal unipotent subgroup of the \fd\ group $G(\bc [z]/z^k)$. This
  inverse system is compatible with the $\bc ^\times$ action. Moreover
  the exponential map $\frak{u}[k] \lra U[k]$ is a $\bc
  ^\times$-equivariant \iso\ of varieties, so we may identify $U[k]$
  with a \fd\ representation of $\bc ^\times$. 

The next proposition is a special case of \cite{kumar}, Proposition
7.3.7.

\begin{proposition} \label{} \marginpar{}

Fix $m$. Then for all $k>>0$, the natural map $F_m \util\lra
U[k]$ is a closed \tchat -equivariant embedding of varieties. 

\end{proposition}

We remark that the proof in our special case is quite easy. Using an
embedding $\gc \subset SL_n \bc$ for some $n$, one first reduces to the
case $\gc =SL_n \bc$. Then it is clearly sufficient to take $k>m$.

\bigskip

\begin{corollary} \label{} \marginpar{}

$F_m \util$ is $\bc ^\times$-equivariantly
  isomorphic to a positively weighted cone.
Moreover the inclusions $F_m \subset F_{m+1}$ are induced by inclusions
  of $\bc ^\times$ representations. 
\end{corollary}

\begin{corollary} \label{} \marginpar{}

Let $\tilde{U}_\Theta $ be the ind-subgroup of \util\ associated to a
bracket closed subset $\Theta$ of the positive affine roots
(cf. \cite{kumar}, \S 6.1.1). Then $\tilde{U}_\Theta$ has a commensurate
filtration by irreducible positively weighted cones.

\end{corollary}

\proof It is clear that $\tilde{U}_\Theta$ is a connected ind-group,
hence an irreducible ind-variety by \cite{kumar}, Lemma 4.2.5. Since
$\tilde{U}_\Theta$ is $\bc ^\times$-invariant, it inherits a filtration
by weighted cones and hence a commensurate filtration by irreducible
weighted cones. 

\begin{corollary} \label{attractpoint} \marginpar{}

The identity element $1$ is strongly attractive in \util .

\end{corollary}

\proof Consider the case \util , and choose an open $V \subset \util$
so that each $V_m$ is a \nbhd\ of $1$ in $F_m \util$ with compact
closure. Then for any \nbhd\ \cW\ of $1$, $V$ flows to $\cW$ at
zero. Since every $x \in \util$ lies in such a $V$, this proves the
corollary. The case \utilminus\ is the same. 

\subsubsection{Joins}                                                           
Let $X$ be a positively weighted cone in $V$. Choose a Hermitian metric
invariant under the $S^1$ action, and let $S(V)$, $D(V)$ denote
respectively the unit sphere and unit disc. Let $S(X)=X \cap S(V)$,
$D(X)=X \cap D(V)$. Then it is clear that the map $(S(X) \times
[0,\infty))/(S(X) \times \{0\}) \lra X$ given by $(v,t) \mapsto t \cdot
v$ (if $t>0$) and $(v,0) \mapsto 0$ is a $\bc ^\times$-equivariant
\homeo . In particular, $S(X)$ is an equivariant deformation retract of
$X -0$, and $D(X)$ is just the cone $CS(X)$ on $S(X)$---in the
topologist's sense, where $CY =(Y \times [0,1])/(Y \times 0)$.

Now recall that the {\it join} of spaces $Y,Z$ is defined by $Y*Z =(CY
  \times Z) \cup _{Y \times Z} (Y \times CZ)$. An elementary argument
  shows that $Y *Z$ is a \defr\ of $(CY \times CZ)-(p,q)$,
  where $p,q$ are the cone points. Here we conclude:

\begin{lemma} \label{join} \marginpar{}                                         
Suppose $X \subset V$, $Y \subset W$ are positively weighted cones. Then
$(X \times Y)- (0,0)$ contains $S(X)*S(Y)$ as a \defr .

\end{lemma}                                                                     
\subsection{Punctured Birkhoff strata}

The flow shows immediately that a Birkhoff stratum \slam\ itself is
contractible. {\it A priori}, however, there are no restrictions
whatever on the homotopy type of a contractible space minus a point; one
has only to think of the cone on a space minus the cone point. However: 

\begin{lemma} \label{puncture} \marginpar{}

Every punctured Birkhoff stratum $\slam -\{\lam \}$ is   
contractible.

\end{lemma}  

Since $\slam -\{\lambda\}$ is isomorphic as an ind-variety to
\utilaminus , with $\lambda$ corresponding to the identity $e$, it will
be enough to prove:

\begin{lemma} \label{} \marginpar{}

For all $\lambda \in \wtili $, $\utilamp -e$ and $\utilaminus - e$ are
contractible.

\end{lemma}

\proof We consider the case $\utilamp -e$. Suppose for
convenience that $I=\emptyset$, so that $\wtili =\wtil$. Let $\mu \in
\wtil$ satisfy $\el (\lambda \mu) =\el (\lambda) + \el (\mu)$. Then by
general results from \cite{kumar} (see especially Theorem 5.2.3c) group
multiplication defines an \iso\ of ind-varieties

$$\phi : \lambda \tilde{U} _\mu \lambda ^{-1} \times \tilde{U}^\prime _{\lambda
  \mu} \stackrel{\cong}{\lra} \utilamp .$$

\noindent Note that $\phi$ can also be interpreted in terms of {\it
loc. cit.}, Lemma 6.1.3, writing $\phitilplus$ as the disjoint union of
suitable bracket closed subsets.

Now $\tilde{U} _\mu $ is a \fd\ unipotent group of dimension $\el
(\mu)$, where $\el (\mu)$ can be taken arbitrarily large. Moreover,
there is an analogous isomorphism for general $I$. We conclude that for
every $n>0$ there exists $d \geq n$ such that \utilamp\ has a
commensurate filtration by varieties of the form $\bc ^d \times A_n$,
where $\bc ^d$ and $A_n$ are positively weighted cones. Hence by
Lemma~\ref{join}, the corresponding filtrations of $\utilamp -
e$ have the homotopy type of $S^{2d-1} * S(A_n)$, which in turn is
homotopy equivalent to $S^{2d} \wedge S(A_n)$ (\cite{hatcher}, Exercise
0.24) and hence is $(2d-1)$-connected.  Passing to the direct limit,
$\utilamp -e$ is $(2d-1)$-connected. But $d$ can be taken
arbitrarily large, and therefore $\utilamp -e$ is weakly
contractible. Since $\utilamp -e$ is an ind-variety and hence an
ind-CW-complex, this completes the proof.

\section{Schubert and Birkhoff \nbhd s}

Let \cI\ be an upper order ideal and \cJ\ a lower order ideal for the
Bruhat order on \wtili . Then the {\it Schubert neighborhood} $\vi =\cup
_{\lambda \in {\cal I}} \elam$ is a Zariski open \nbhd\ of \zi , and the
{\it Birkhoff neighborhood} $\sj =\cup _{\lambda \in {\cal J}} \slam$ is
a Zariski open \nbhd\ of \xj\ (see the appendix). Recall from \S 2 that
\zi\ is a finite union of \bvars , while \xj\ is a finite union of
\svars . Although we are mainly interested in the case of principal
order ideals---i.e., in Birkhoff and Schubert varieties---the general
case will be useful for later induction arguments.

\begin{theorem} \label{tube} \marginpar{}

a) Let \zi\ be a finite union of \bvars . Then \zi\ is a \defr\ of its
Schubert \nbhd\ \vi .

b) Let \xj\ be a finite union of \svars . Then \xj\ is a \defr\ of its
Birkhoff \nbhd\ \sj . 

\end{theorem}

\proof For ease of notation, we write $Z, \cE, X, \cS$ in place of \zi ,
\vi , \xj , \sj .

a) By Theorem~\ref{flow}, it is enough to show that $Z$ is
strongly attractive in \cE , or that for every \svar\ $X=X_\lambda$, $Z
\cap X$ is strongly attractive in $\cE \cap X$. Let $f: \tilde{U} \times
(Z \cap X) \lra \cE \cap X$ denote the map induced by the action of
\util\ on \cE . Note that $f$ is a $\bc ^\times $-equivariant map of
pairs $(\util \times (Z \cap X), \{1\} \times (Z \cap X)) \lra (\cE \cap
X , Z \cap X)$, where $\bc ^\times$ acts on $\util \times Z$ by $t \cdot
(u,z) =(tut^{-1}, t\cdot z)$. We will prove (a) by showing that $f$
satisfies the hypotheses of Lemma~\ref{star}. 

Since $Z \cap X$ is compact, tubes of the form $V \times (Z \cap X)$ are
cofinal among \nbhd s of $Z \cap X$ in $\util \times (Z \cap X)$. It
then follows from Corollary~\ref{attractpoint} that $\{1\} \times (Z
\cap X)$ is strongly attractive in $\util \times (Z \cap X)$.  Next we
show that $f$ satisfies condition (*) of Lemma~\ref{star}. For each
$\sigma \in I \cap J_\lambda$, the natural map $\tilde{U} _\sigma \times
S_\sigma \lra \cU _\sigma$ is an isomorphism of ind-varieties, and
restricts to an \iso\ $\tilde{U} _\sigma \times (S_\sigma \cap X) \cong
\cU _\sigma \cap X$ (see the Appendix). More generally, for any $g \in
\util$ we have $g\tilde{U} _\sigma \times (S_\sigma \cap X) \cong g\cU
_\sigma \cap X$. Since the ind-varieties $g \tilde{U} _\sigma \times
S_\sigma$ cover $\util \times Z$, this verifies condition (*).

b) The proof here is analogous to the proof of (a), using the flow at
infinity. In this case we use the natural map $f: \utilminus \times X
\lra \cS$. Since $\utilminus$ is an ascending union of negatively
weighted cones, and $X$ is compact, we conclude as before that $\{1\}
\times X$ is strongly attractive at infinity in $\utilminus \times
X$. Condition (*) of Lemma~\ref{star} is also verified as in (a), using
the \iso s $\tilde{U}^- _\sigma \times e_\sigma \stackrel{\cong}{\lra}
\cU _\sigma$. 

\bigskip

Variants of Theorem~\ref{tube} can be obtained by intersecting with
$\bc ^\times$-invariant closed ind-subvarieties of \cF . In particular,
we will need the following for the proof of Theorem~\ref{main}: 

\begin{theorem} \label{tube2} \marginpar{}

If $\cI \subset \cI ^\prime$, then $Z_{{\cal I}} $ is a \defr\ of 
$\cE _{{\cal I}} \cap Z_{{\cal I}^\prime}$. 

\end{theorem}
 
\proof Since $Z_{{\cal I}^\prime}$ is invariant under the flow, the
proof of Theorem~\ref{tube} shows that $Z_{{\cal I}} $ is strongly
attractive in $\cE _{{\cal I}} \cap Z_{{\cal I}^\prime}$.

\section{The homotopy-type of a Birkhoff variety}

\subsection{The main theorem}

In this section we show that every finite union of \bvars\ \zi\ is a
\defr\ of \cF\ (Theorem~\ref{equivalence}). To motivate this result, we
point out that there are much simpler examples of the same
phenomenon. For instance, the ind-variety \cpinf\ has a Birkhoff
filtration $\cpinf \supset Z_1 \supset Z_2 \supset...$ dual to its
Schubert filtration $\bp ^n$: Writing $\bc ^\infty =\cup \bc^n$ as
usual, $Z_n$ is just the subvariety of lines orthogonal to $\bc ^n$, and
hence is isomorphic to \cpinf . Moreover it is easy to show that
$Z_{n+1}$ is a \defr\ of $Z_n$, as follows: The ind-variety structure on
$Z_n$ is given by the {\it Richardson varieties} $Z_n \cap \bp ^m$. But
$Z_n \cap \bp ^m$ is just $\bp ^{m-n}$, and hence $(Z_n \cap \bp
^m)/(Z_{n+1} \cap \bp ^m) =S^{2(m-n)}$. Letting $m \ra \infty$, we have
$Z_n /Z_{n+1} =S^\infty$, which is contractible. Since $(Z_n, Z_{n+1})$
is a CW-pair, it follows that $Z_{n+1}$ is a \defr\ of $Z_n$. In
principle one could apply the same method in the present context, but
the Richardson varieties $\zlam \cap X_\sigma$ are much more
complicated. Hence we will take a somewhat different approach.

\begin{theorem} \label{equivalence} \marginpar{}

Let $Z \subset Z^\prime$ be finite unions of \bvars . Then $Z$ is a \defr\ of
$Z^\prime$. In particular, any \bvar\ \zlam\ is a \defr\ of \cF .

\end{theorem}

\proof By downward induction over the strata, we reduce at once to the
case when $Z^\prime -Z$ is a single stratum \slam , and $Z^\prime$ is a
\defr\ of \cF . Since \cF\ is simply-connected (it is a connected
CW-complex with only even-dimensional cells), in particular $Z^\prime$
is simply-connected. Let \cE\ be the Schubert \nbhd\ of $Z$. Then we
have a diagram of open sets and inclusions 

\bigskip                                                                        
                                                                                
$                                                                               
\begin{diagram}                                                                 
\node{\slam -\{\lam \}} \arrow{e,t}{i} \arrow {s,l}{}                            
\node{\slam} \arrow{s,r}{} \\                                                   
\node{\cE \cap Z^\prime} \arrow{e,b}{j} \node{Z^\prime}                          
\end{diagram}                                                                   
$                                                                               
                                                                                
\bigskip         

\noindent where \slam , $\cE \cap Z^\prime$ cover $Z^\prime$ and have
intersection $\slam -\{\lambda\}$. Then $H_*i$ is an \iso\ by
Lemma~\ref{puncture}, so $H_* j$ is an \iso\ by excision. Since $\slam
-\{\lambda \}$, \slam\ and $\cE \cap Z^\prime$ are path-connected (for
the last case see Proposition~\ref{pathconnect}), we conclude similarly
from the Seifert-van Kampen theorem that $\pi _1 j$ is an \iso\ and
hence $\cE \cap Z^\prime$ is simply-connected. But $Z$ is a \defr\ of
$\cE \cap Z^\prime$ by Theorem~\ref{tube}a. Hence $Z$ is a \defr\ of
$Z^\prime$ by Corollary~\ref{wedefr}.

\bigskip

It follows, of course, that the inclusions $Z \subset \cF$ induce \iso
s on any homology or cohomology theory, including \that -equivariant
theories. For emphasis we record the following cases explicitly.

\begin{corollary} \label{} \marginpar{}

For any finite union of \bvars\ $Z$, 
$H^* \cF \stackrel{\cong}{\lra} H^* Z$, and 
$H^* _{\hat{T}} \cF \stackrel{\cong}{\lra} H^* _{\hat{T}} Z$.

\end{corollary}

\noindent {\it Remark:} It follows from the corollary that $H_* Z$ has
finite type and is concentrated in even dimensions. This does not seem
obvious {\it a priori}; neither property need hold for an
ind-subvariety $Y$ of \cF . 

\subsection{Cohomology of $(\sj , \sj - \xj )$ and $(\vi , \vi -\zi )$}

We next consider the pairs $(\sj , \sj - \xj )$ and $(\vi , \vi -\zi )$,
which can be viewed as ``normal Thom spaces'' of the subvarieties $\xj ,
\zi$ in \cF .  Let \cI\ be a nonempty upper order ideal and let \cJ\
be the complementary lower order ideal. Then

$$\cF =\cE _{{\cal I}} \cup \sj ,$$ 

$$\vi - \zi =\vi \cap \sj =\sj -\xj .$$

For ease of notation, we henceforth write $\cE , \cS , Z, X$ for the
corresponding spaces above. 

\begin{proposition} \label{} \marginpar{}

$H_* (\cS , \cS -X) =0$. Hence $H_* (\cS -X) \cong H_* X$. 

\end{proposition}

\proof We have 

$$H_* (\cS , \cS -X) \cong H_* (\cF , \cE )\cong H_* (\cF , Z)=0,$$

\noindent where the first \iso\ is by excision, the second by
Theorem~\ref{tube2} and the third by Theorem~\ref{equivalence}. Thus
$H_* (\cS -X) \cong H_* \cS \cong H_*X$ by Theorem~\ref{tube2}. 

\bigskip

\noindent {\it Remark:} This result reflects the intuition that \cS\ is
a sort of infinite-dimensional ``vector bundle'' over $X$, and so should have
contractible Thom space, while its ``sphere bundle'' should have contractible
fibers. 

\bigskip

Similarly, we have: 

\begin{proposition} \label{} \marginpar{}

$H_* (\cE , \cE -Z) \cong H_* (\cF , X)$, and $H_* (\cE -Z) \cong H_*
  X$. 

\end{proposition}

\proof We have 

$$H_* (\cE , \cE -Z) \cong H_* (\cF , \cS) \cong H_* (\cF , X).$$

\bigskip

\noindent {\it Remark:} Note that $H_* (\cF , X)$ is a free abelian
group on the upper order ideal \cI , graded by twice the length as
usual. In fact $\cF /X$ is a CW-complex whose cells are the Schubert
cells corresponding to \cI , plus a basepoint. When $\cI =\cI _\lambda$,
$\cF /X$ has $2 \el ^S (\lambda)$-skeleton the sphere $e_\lambda
^+$. This reflects the intuition that the pair $(\vi , \vi - Z
_\lambda)$ is the ``Thom space'' of the complex $\el ^I (\lam
)$-dimensional ``normal bundle'' of \zlam\ in \affg . In cannot actually
be such a Thom space, however, since it does not even have the right
\pser .

\subsection{Remarks on equivariant cohomology}

Let $Y$ be a space with an action of a compact torus $T$. The
Goresky-Kottwitz-MacPherson (GKM) theory \cite{gkm} characterizes the
equivariant cohomology $H_T^* (Y; \bq)$ in terms of the zero- and
one-dimensional orbits---provided that $Y$ is sufficiently well-behaved
as a $T$-space. In particular, some finiteness restriction on $Y$ is
usually required, such as compactness, finite cohomological dimension,
and/or finite orbit type. Since the spaces we are considering are
noncompact, of infinite cohomological dimension, and of infinite orbit
type, any attempt to extend the results of \cite{gkm} must proceed with
caution.

The case of the \that -action on \cF\ itself has been studied by a
number of authors; see \cite{kumar} and the references cited there, and
\cite{hhh}. Here the beautiful properties of the Schubert cell
decomposition more than compensate for the infinite-dimensionality of
\cF ; one can proceed by induction over the Schubert filtration. The
result is as follows: Identify $H^* _{\hat{T}} (\cF ^{\hat{T}})$ with
the ring of $H^*_{\hat{T}}$-valued functions on $\wtili $. Let $\cR
(\cF )$ denote the subring consisting of those functions $f$ such that
whenever $r_\theta \sigma =\lambda$ for some positive affine root
$\theta$, we have $f(\sigma ) =f(\lambda ) \, mod \, c_\theta$, where
$c_\theta$ is the first Chern class of the line bundle $\xi _\theta
\downarrow B\that$ associated to $\theta$. Then \cF\ satisfies the
{\it GKM theorem} (compare \cite{gkm}, 1.2.2); that is, restriction to
the fixed point set defines an \iso\ $H^* _{\hat{T}} (\cF ;\bq) \cong
\cR (\cF ; \bq)$.

Now consider a Birkhoff subspace $Z=\zi$. Again we are faced with an
infinite-dimensional space, with the further complication that there is
no Schubert cell structure. The natural filtration by Richardson
varieties is not so easy to analyze. Instead we will use
Theorem~\ref{equivalence} to obtain half of the GKM theorem for $Z$.
Note that $Z^{\hat{T}}=\cI$, and that if $\sigma > \lambda \in \cI$,
then the unique one-dimensional orbit with $\sigma, \lambda$ as its
poles lies in $\zlam \cap X_\sigma \subset Z$.

\begin{proposition} \label{equi} \marginpar{}

For any Birkhoff subspace $Z=\zi$, the restriction $i^*: H^*_{\hat{T}} Z \lra
H^* _{\hat{T}} Z^{\hat{T}}$ is an injection into $\cR (Z)$. 

\end{proposition}

\proof That $i^*$ has image contained in $\cR (Z)$ is straightforward;
the argument is as in \cite{gkm} or \cite{kumar}. It remains to show
that $i^*$ is injective. Since $H^* _{\hat{T}} \cF$ is torsion-free,
it suffices to prove this rationally. Let $j: Z \lra \cF$ denote the
inclusion. Then there is a commutative diagram

\bigskip

$
\begin{diagram}
\node{H^*_{\hat{T}} (\cF ;\bq)} \arrow{e,t}{\cong} \arrow {s,l}{\cong}
\node{H^*_{\hat{T}} (Z;\bq)} \arrow{s,r}{i^*} \\
\node{\cR (\cF ;\bq)} \arrow{e,b}{\cR (j)} \node{\cR (Z;\bq)}
\end{diagram}
$

\bigskip

\noindent Thus $i^*$ is injective if and only if $\cR (j)$ is
injective. Now $\cR (j)$ amounts to taking a function $f: \wtil/\wtilil
\lra H^* _{\hat{T}}$ and restricting it to the upper order ideal \cI
. Suppose that $f$ restricts to zero on \cI , and let $\sigma \wtilil
\in (\wtil /\wtilil) -\cI$. Let \cA\ denote the set of all affine
reflections. Since \cA\ is infinite and all isotropy groups of the
action of \wtil\ on $\wtil /\wtilil$ are finite, the set $(\cA 
\sigma \wtilil ) \cap \cI$ is infinite. But then $f(\sigma \wtilil)$ is
divisible by an infinite set of pairwise relatively prime elements of
$H^2 (B\hat{T}; \bq)$, and hence must be zero. This proves the
proposition.

\bigskip

From the commutative diagram we also have $Im \, i^* =Im \, \cR (j)
$. Hence the full GKM theorem holds if and only if every function $f \in
\cR (Z; \bq)$ extends to $\tilde{f} \in \cR (\cF ; \bq)$. 

\bigskip

\noindent {\it Remark:} The \that -space $\cE \cap \cS$ provides a
typical example of what can go wrong with equivariant cohomology in an
infinite-dimensional setting. It is equivariantly formal in the sense of
\cite{gkm}, since $H^* _{\hat{T}} (\cE \cap \cS)$ is a free module $H^*
_{\hat{T}} \otimes H^* X$, but it has no $\hat{T}$-fixed points.  Hence
localization at the fixed point set and the GKM theorem fail for $\cE
\cap \cS$.

\section{Appendix: Basic properties of \bvars }

We assume given the standard refined Tits system structure on $\tilde{G}$;
in particular, the Bruhat and Birkhoff decompositions (\cite{kp},
\cite{kumar}). Recall that $\openlam =\lam \cU _0$ (\S 2), where $\cU _0
=\tilde{U}_I ^- P_I/P_I$. 

\begin{proposition} \label{} \marginpar{}

The \openlam 's form a Zariski open cover of $\cF $. 

\end{proposition}

\proof That the \openlam 's cover $\cF$ follows from the Bruhat
decomposition. By reducing to the case $G=SL_n \bc$, it is not hard to
show that $\cU _0$ is Zariski open (or see \cite{kumar}). Since
multiplication by any fixed $f \in \tilde{G}$ gives a morphism of
ind-varieties from \cF\ to itself, it follows that all the \openlam 's
are Zariski open.

\bigskip

\begin{proposition} \label{product} \marginpar{}

The natural maps $\utilam \times \slam \lra \openlam$ and $\utilaminus
\times \elam \lra \openlam$ are \iso s of ind-varieties. 

\end{proposition}

\noindent {\it Proof:} The first \iso\ amounts to the \iso\ of
ind-varieties 

$$\phi: (\util \cap (\lam \tilde{U}_I ^- \lambda ^{-1})) \times
(\utilminus \cap (\lam \tilde{U}_I ^- \lambda ^{-1})) \cong \lam
\tilde{U}_I ^- \lambda ^{-1}.$$

\noindent Here $\phi$ is group multiplication. That $\phi$ is bijective
follows from the axioms for a refined Tits system; compare \cite{kumar},
p. 169 (7), as well as p. 227 (1). The methods there also show that
$\phi$ is an \iso\ of ind-varieties. 

\begin{corollary} \label{intcor} \marginpar{}

$U_\lambda \times (\slam \cap \xmu ) \lra \openlam \cap \xmu$ is an
  \iso\ of varieties.

\end{corollary}

\bigskip

If \cJ\ is a lower order ideal in $\wtils /W$, let $\sj =\cup _{\lambda
\in {\cal J}} \slam $. If \cI\ is an upper order ideal, let $\vi =\cup
_{\lambda \in {\cal J}} \elam$. And if \cK\ is any subset of $\wtils /W$, let
$\cU _{{\cal K}} =\cup _{\lambda \in {\cal K}} \openlam $.

\begin{proposition} \label{snbhd} \marginpar{}

a) Let \cJ\ be a lower order ideal. Then $\xj \subset \sj =\cU
_{{\cal J}}$. 

b) Let \cI\ be a upper order ideal. Then $\zi \subset \vi =\cU
_{{\cal I}}$. 

\end{proposition}

\proof In (a) we have $\xj \subset \cU _{{\cal J}}$ and $\sj\subset \cU
_{{\cal J}}$ by Proposition~\ref{product}. Now suppose $\lam \in \cJ$;
we show that $\cU _{\lambda} \subset \sj $. Since $\sj$ is \utilminus
-invariant, it is enough to show $\elam \subset \sj $. But if $x \in
\elam \cap S_\mu$, then $\mu =lim _{t\ra \infty} t\cdot x \in
\overline{e} _\lambda$, so $\mu \leq \lam $.

The proof of (b) is similar. 

\bigskip

\begin{corollary} \label{} \marginpar{}

Let \cI , \cJ\ be respectively upper and lower order ideals. Then

a) $\sj$ is Zariski open and $\zi$ is Zariski closed. 

b) $\vi$ is Zariski open and $\xj$ is Zariski closed. 

\end{corollary}

Both statements follow immediately, using the fact that the complement
of an upper order ideal is a lower order ideal and vice-versa. 

\bigskip

\begin{proposition} \label{closure} \marginpar{}

$\overline{S}_\lambda =\zlam$ and $\overline{e} _\lambda =\xlam $. 
Here the closure can be taken in either the Zariski topology or the
classical topology. 

\end{proposition}

\noindent {\it Proof:} By the corollary, we have $\overline{S}_\lambda
\subset \zlam$ and $\overline{e} _\lambda \subset \xlam $. The reverse
inclusions reduce to showing that if $\sigma \da \eta$, then $\sigma \in
\overline{S}_\eta$ and $\eta \in \overline{e} _\sigma$. Let $r_\theta
\sigma =\eta$ and let $SL_2 ^\theta \subset \tilde{G}$ denote the
corresponding $SL_2$ subgroup. Then $SL_2 ^\theta \cdot \sigma$ is an
embedded $\bp ^1$, denoted $\bp ^1 _{\eta \sigma}$, with $\sigma , \eta
\in \bp ^1$ and $\bp ^1 -\{\sigma, \eta \} \subset S_\sigma \cap
e_\eta$. This proves the proposition.

\bigskip

\begin{proposition} \label{} \marginpar{}

If $\lam \leq \mu$, then $\zlam \cap \xmu$ is irreducible of codimension
$\el ^I \lam$ in \xmu . 

\end{proposition}

\proof Note that $\openlam \cap \xmu$ is irreducible, since it is
Zariski open in the irreducible variety \xmu . Hence $\slam \cap \xmu$
is irreducible of codimension $\el ^I \lam$ by
Corollary~\ref{intcor}. Since $\slam \cap \xmu$ is Zariski dense in
$\zlam \cap \xmu$, the result follows. 

\begin{proposition} \label{pathconnect} \marginpar{}                                
Suppose \cI , \cJ\ are respectively upper and lower order ideals in
\wtili . If $\cI \cap \cJ$ is connected as a subgraph of the Hasse
diagram of \wtili , then $\zi \cap \xj $ is path-connected. In
particular $\zi \cap \xlam$, $\zlam \cap \xj$, $\zi$, and $\xj$ are all
path-connected.

\end{proposition}                                                               
\proof Clearly the strata and cells are path-connected. If $\cI \cap \cJ$ is
connected in the Hasse diagram, then any two of its points can be joined
by a sequence of $\bp ^1 _{\lambda \mu}$'s (see the proof of
Proposition~\ref{closure}) lying in $\zi \cap \xj$.

\end{document}